\documentclass[10pt,twoside]{article}

\usepackage{amsmath,amsthm,verbatim,amssymb,amsfonts,amscd,diagrams,graphics}

\oddsidemargin=\evensidemargin
\newtheorem{theorem}{Theorem}[section]
\newtheorem{corollary}[theorem]{Corollary}

\newtheorem{proposition}[theorem]{Proposition}

\theoremstyle{definition}

\theoremstyle{remark}
\newtheorem{remark}[theorem]{Remark}

\newcommand{\td}[1]{\tilde{#1}}
\newcommand{\into}{\hookrightarrow}
\newcommand{\Z}{\mathbb{Z}}

\newcommand{\R}{\mathbb{R}}
\newcommand{\bd}{\partial}

\newcommand{\bG}{\bar G}
\newcommand{\bg}{\bar g}

\newcommand{\mc}[1]{\mathcal{#1}}

\newcommand{\Hom}{\operatorname{Hom}}
\newcommand{\ext}{\operatorname{Ext}}

\newcommand{\mf}{\mathfrak}
\newcommand{\beq}{\begin{equation*}}
\newcommand{\eeq}{\end{equation*}}

\begin{document}
\title{Groups of locally-flat disk knots and non-locally-flat sphere knots}
\author{Greg Friedman\\
 Yale University, Dept. of Mathematics\\
 10 Hillhouse Ave., PO Box 208283. New Haven, CT 06520\\friedman@math.yale.edu - Tel. 203-432-6473  Fax:  203-432-7316}
%11/22/02
\maketitle
%typeset=\today
%\vskip.75in
\begin{abstract}

The classical knot groups are the fundamental groups of the complements of smooth or piecewise-linear (PL) locally-flat knots. For PL knots that are not locally-flat,  there is a pair of interesting groups to study: the fundamental group of the knot complement and that of the complement of the ``boundary knot'' that occurs around the singular set, the set of points at which the embedding is not locally-flat. If a knot has only point singularities, this is equivalent to studying the groups of a PL locally-flat disk knot and its boundary sphere knot; in this case, we obtain a complete classification of all such group pairs in dimension $\geq 6$. For more general knots, we also obtain complete classifications of these  group pairs under certain restrictions on the singularities. Finally, we use spinning constructions to realize further examples of boundary knot groups. 

\end{abstract}

\section{Introduction}

In the author's dissertation (see \cite{GBF} and \cite{GBF1}), we studied the generalization of Alexander polynomials to PL sphere knots which were not necessarily locally-flat, i.e. PL-embedding $S^{n-2}\into S^n$ such that the neighborhood disk pairs of points in the image of the embedding are not necessarily PL-homeomorphic to the standard disk pair. In this paper, we study the generalization to such knots of another classical knot invariant, the knot group. 

The classically studied knot groups are the fundamental groups of the complements of smooth or PL locally-flat knots. In the context of our PL singular knots, there is a pair of interesting groups to study: the fundamental group of the knot complement and that of the complement of the ``boundary knot'' which occurs around the \emph{singular set}, the set of points at which the embedding fails to be locally flat. If a knot has only point singularities, this is equivalent to studying the group of a locally-flat disk knot and that of its boundary locally-flat sphere knot, and in this case, we obtain a complete classification of all such group pairs for knots of dimension $n\geq 6$. For more general knots, we also obtain complete classifications of the main knot group and of the boundary knot group under certain restrictions on the singularities. Finally, we show how spinning constructions can be used to realize further examples of boundary knot groups. Note that all embeddings in this paper are Piecewise Linear (PL).

We now outline our results in slightly greater detail:

The groups of smooth or PL locally-flat sphere knots $K: S^{n-2}\into S^n$ were completely classified for $n\geq 5$ by Kervaire in \cite{Ke}. In Section \ref{S: groups} of this paper, we show in Theorem \ref{T: nec con} that Kervaire's necessary conditions extend to all PL-knots. In Theorem \ref{T: class} and Corollary \ref{C: decouple}, we obtain a classification, analogous to Kervaire's,  for the pair of groups associated to a locally-flat disk knot $J: D^{n-2}\into D^n$, $n\geq 6$ and its boundary locally-flat sphere knot. This implies a classification in the same dimensions of groups for sphere knots with point singularities. 

 In Section \ref{S: links}, we study the \emph{boundary knots} of singular sphere knots with singular sets of dimension $>0$. These boundary knots are not necessarily sphere knots, but they will be locally-flat codimension-2 manifold pairs. In particular, we show that the boundary knot groups will not, in general, satisfy the Kervaire conditions, but in Theorem \ref{T: bd kerv}, we establish that they will if the singular set is $2$-connected. We also compute the homology of  these boundary knot complements in terms of the homology of the singular set (Theorem \ref{T: H of X}). In case the singular set  has a single stratum, we establish some further necessary conditions on the boundary knot group in terms of the fundamental groups of the stratum and its ``link knot'' (Theorem \ref{T: manif sing} and its corollaries).

Finally, in Section \ref{S: realize}, we show how to realize some further examples of knot group pairs via knot constructions such as frame twist-spinning and suspension.

\section{Preliminaries and conventions}\label{S: cons}

\paragraph{Basic definitions} We define a knot, $K$, to be a PL-embedding $S^{n-2}\into S^n$; we do not assume that the embedding is locally-flat, i.e. there may be points whose regular neighborhoods pairs are not PL homeomorphic to the standard unknotted disk pair. Following standard abuse of notation, we sometimes also use $K$ to refer either to the image $K(S^{n-2})$ or the pair $(S^n, K(S^{n-2}))$. We also sometimes simply refer to the knotted sphere pair $(S^n, K)$. We let $\Sigma$ denote the \emph{singular set} of points at which the embedding $K$ fails to be locally-flat, and without further mention we identify $\Sigma$ either as a subset of $S^{n-2}$ or of $K(S^{n-2})\subset S^n$. 

Given a knot, we will most often be concerned with the topological properties of its complement $S^n-K$, which is homotopy equivalent to the knot exterior,  $C$, which is the complement in $S^n$ of an open regular neighborhood $N(K)$ of $K$ (generally, we let $\bar{N}$ stand for closed regular neighborhoods and $N$ stand for the interior of $\bar N$). In case the embedding is not locally-flat, we will also be concerned with the $\emph{boundary knot}$. If $\Sigma$ is the singular set of the embedding,  then the boundary knot is the pair $(\bd \bar N(\Sigma), \bd \bar N(\Sigma)\cap K)$, where $\bar N(\Sigma)$ is the closed regular neighborhood of $\Sigma$ in $S^n$. Note that the boundary knot does not necessarily consist of knotted spheres, but it is a locally-flat codimension two manifold pair. We can then consider the boundary knot complement $\bd \bar N(\Sigma)-(\bd \bar N(\Sigma)\cap K)$ and the homotopy equivalent exterior $X$, the complement in $\bd \bar N(\Sigma)$ of an open regular neighborhood of $\bd \bar N(\Sigma)\cap K$. It is not hard to see that we can choose these regular neighborhoods such that $X$ is a subspace of $C$. In fact, $C$ and $X$ are manifolds with boundary, and $\bd C\cong X\cup_{\bd X=\bd T} T$, where $T$ is a circle bundle over the manifold $K-( N(\Sigma)\cap K)$ (see below for more about this property of $T$). 
Our main objects of study will be the fundamental groups of $C$ and $X$ and the homomorphism between them induced by the inclusion $X\into C$. 

\paragraph{The relationship between knots with point singularities and disk knots} In the special case where $K$ is a knot with a single point singularity, then choosing  $N(\Sigma)$ as the star of $\Sigma$ in the second barycentric subdivision of the triangulation, we can identify $S^n-N(\Sigma)$ as an $n$-disk, and the pair $(S^n-N(\Sigma), K\cap (S^n-N(\Sigma))$ as a locally-flat PL disk knot, i.e. a proper locally-flat PL embedding $J: D^{n-2}\into D^n$. In this case, the boundary knot is simply the boundary PL locally-flat sphere knot of the disk knot. If $K$ has multiple point singularities, we can use a technique of Fox and Milnor \cite{M66} to slightly modify the definition of the boundary knot to obtain  again a nice disk knot pair: Let $\rho$ be a simplicial path with no crossings in $K(S^{n-2})$ that connects the singular points of the knot, in other words, a path that starts at one singular point and then traverses all of them in some order with no self-intersections. Instead of the closed regular neighborhood $\bar N(\Sigma)$, we can instead consider the regular neighborhood $\bar N(\rho)$. Since $\rho$ is contractible, this neighborhood is a disk, and the complement of its interior again gives a disk knot. In this case, the boundary knot is the knot sum of the link knots about the point singularities. Note that for our purposes this construction is essentially independent of the choice of $\rho$: Since the regular neighborhood of $\rho$ collapses into the knot $K$, we can see that the disk knot exterior we obtain by this construction is isomorphic to the knot exterior $S^n-N(K)$. Similarly, the boundary knot complement $X$ will be the complement of the sphere knot given by the knot sum of the link knots around the singular points, and this depends only on the knot sum itself, not on the order in which we connect these knots, because knot sum is associative and commutative for locally-flat sphere knots.
\begin{comment}
 The map on $\pi_1$ determined by inclusion $X\to C$ will also be unaffected by choice of $\rho$ because $\pi_1$ of the complement of a knot sum will be a free product with amalgamation over $\Z$ of the knot groups of the summands, and the map is therefore determined by its action on the generators of $\pi_1$ of each knot summand. But the map on $\pi_1$ determined by the inclusion of a path representing an element of $\pi_1$ of a summand into $C$ is clearly independent of $\rho$. (Of course this whole discussion is predicated upon choosing some fixed set of basepoints in all of the relevant spaces and appropriate reference paths between them; the standard  ambiguity can indeed creep in here, but, as 
usual, only up to certain isomorphisms).
\end{comment}

If we are given a locally-flat disk knot, there is a converse to the above construction which will give us a sphere knot with point singularity: we can simply add the cone pair on the boundary. However, this construction only gives us knots with a single point singularity. A more general construction would be the following: If a locally-flat disk knot $J: D^{n-2}\into D^n$ has a  boundary locally-flat sphere knot $K: S^{n-3}\into S^{n-1}$ which can be written as  a connected sum of knots $K=K_1\#\cdots \# K_k$, then there is an ambient isomorphism of the boundary $S^{n-1}$ which can be extended to the interior of the disk and which  arranges the knot so that each $K_i$ is contained in a disk except for the tubes which connect it to the other knots in the sequence. In other words,  we can assume the sum $K_1\#\cdots \# K_k$ represents the knot sum embedded in the standard way, where we begin with the summand knots completely separated within non-intersecting balls $D^{n-1}_i$ in $S^{n-1}$ and then connect them via non-intersecting and non-self-intersecting tubes $D^1\times S^{n-4}$ (technically, a surgery on the knots determined by $1$-handles embedded in general position in $S^n$, or, equivalently, an internal connected sum). We can assume that the tubes intersect each $\bd D^{n-1}_i$ in general position  and furthermore the tube connecting $K_i$ to $K_{i+1}$ will intersect only $\bd D^{n-1}_i$ and $\bd D^{n-1}_{i+1}$, each only once and each intersection being PL-homeomorphic to $S^{n-4}$. Then we can construct a knot with $k$ singularities whose link knots are the $K_i$ as follows: For the ambient sphere, we simply add the cone on the boundary $S^{n-1}$. For the knotted sphere, we take the union of the disk knot $J$, $k$ separate non-intersecting cones on the sets $K\cap D^{n-1}_i$, and  the connecting tubes, which we fill in outside of the interiors of the $D^{n-1}_i$ to form closed tubes $D^1\times D^{n-3}$. As the basepoint of each cone, we can take an arbitrary point on the open cone line from the center of each $D^{n-1}_i$ to the cone point of the cone we have appended onto $\bd D^{n}$ to create the ambient sphere. This construction gives us a sphere knots with point singularities whose link knots are the $K_i$.  

\begin{comment}See the  schematic diagram Figure \ref{F: cone}. 

\begin{figure}[h]
\begin{center}
\scalebox{0.3}{\includegraphics{cone.ps}}
\end{center}
\caption{Coning off a knot sum}
\label{F: cone}
\end{figure}
\end{comment}

Therefore, by the preceding paragraphs, the study of the homotopy properties of sphere knots with point singularities can be considered equivalent to the study of such properties of locally-flat disk knots. We assume such an identification throughout the following. 

\paragraph{Knots as stratified spaces} We conclude this preliminary section with some definitions from the theory of stratified spaces that will be useful in the second half of the paper. We first provide some general definitions and then demonstrate how they will apply to the study of  knots. Let us begin by recalling the definition of a stratified pair of paracompact Hausdorff spaces
$(Y,Z)$ as given in \cite{CS}. Let $c(A)$ denote the open cone on the space $A$,
and let $c(\emptyset)$ be a point. Then a \emph{stratification} of $(Y,Z)$ is a filtration
\begin{equation*}
Y=Y_n\supset Y_{n-1} \supset Y_{n-2}\supset \cdots \supset Y_0\supset Y_{-1}=\emptyset
\end{equation*}
such that for each point $y\in Y_i-Y_{i-1}$ (if it is non-empty), there exists a \emph{distinguished neighborhood} 
$N$, a compact Hausdorff pair $(G,F)$, a filtration 
\begin{equation*}
G=G_{n-i-1}\supset  \cdots \supset G_0\supset G_{-1}=\emptyset,
\end{equation*}
and a homeomorphism
\begin{equation*}
\phi: \R^i\times c(G,F)\to (N,N\cap Z)
\end{equation*}
that takes $\R^i\times c(G_{j-1},G_{j-1}\cap F)$ onto $(Y_{i+j},Y_{i+j}\cap Z)$.  This condition says that neighborhoods of points are locally cone bundles over euclidean space. The sets $Y_i$ are called the \emph{skeleta} of $Y$ and the sets $Y_i-Y_{i-1}$ are the \emph{strata}. The definition implies that the stratum $Y_i-Y_{i-1}$ is a manifold of dimension $i$. The pair $(G,F)$ occuring in the definition of a distinguished neighborhood is called the \emph{link} or \emph{link pair} of the point $y$. 

For $(Y,Z)$ a
compact PL pair, such a stratification exists with each $\phi$ a PL map and with the filtration refining
the filtration by $k$-skeletons (see \cite{Bo}). Note, however, that the choice of such a stratification is generally very non-unique (for example, given any triangulation of $Y$ for which $Z$ is a subcomplex, one can filter by simplicial skeleta to obtain  a stratification). We refer the reader to \cite[Ch. I]{Bo} for a more comprehensive treatment of PL stratified spaces. 

Now suppose that we have a PL knot $K$ (recall that  by the standard abuse of notation we may use $K$ to stand for the image of the embedding). In this case, we take $Y=S^n$  and $Z=K$, and we can consider stratifications of the pair $(S^n, K)$. One such stratification is obtained by choosing a fixed triangulation of the PL pair $(S^n, K)$ and then letting $S^n_i$, $0\leq i\leq n-2$, be the union  of  the $i$ simplices in $K$.  Note that then $S^n_{n-2}=K$. We also set $S^{n}_{n-1}=K$ and $S^{n}_n=S^n$. This is easily checked to be a stratification, though once again not a unique one. Since $S^n$ and $K\cong S^{n-2}$ are both manifolds, it is not hard to see that for a point $y\in Y_i-Y_{i-1}$, the spaces $G$ and $F$ in the link pair must be PL homeomorphic to spheres of respective dimensions $n-i-1$ and $n-i-3$, and we call this pair the \emph{link knot}. The pair $(G,F)$ may be non-trivially knotted, even non-locally-flatly.  However, for $i\geq n-3$, any embedding $S^{n-i-3}\into S^{n-i-1}$ must be unknotted, and so the neighborhood $N$ will be the standard (open) unknotted ball pair. Thus the set of non-locally flat points must lie in a subcomplex of dimension $\leq n-4$. Letting $\Sigma$ denote the \emph{singular set of non-locally flat points}, we can  thus always find stratifications of the form $S^n_n\supset S^n_{n-1}=S^n_{n-2}=K\supset S^n_{n-3}= S^n_{n-4}=\Sigma \supset S^n_{n-5}\supset \ldots$. N.B. This convention differs slightly from the standard of allowing $\Sigma$ to represent the entire ``singular locus'' $S^{n}_{n-2}$. 

It is possible to continue describing an explicit stratification determined by the minimal dimensions of link knot pairs. In other words, we could let $S^n_{n-4}=\Sigma$ be the set of points whose distinguished neighborhoods can only be described with  link  knots of dimension  $\geq (3,1)$, $S^n_{n-5}$ the set of points whose distinguished neighborhoods can only be described with  link  knots of dimension $\geq (4,2)$, and so on. However, we will not need this kind of refinement, so we omit further details. In the second half of this paper we will be concerned with knots that allow stratifications with certain properties.

\section{Knot groups}\label{S: groups}

With the notation above, we will refer to $\pi_1(C)$ as the \emph{knot group} and $\pi_1(X)$ as the \emph{boundary knot group} of a not necessarily locally-flat knot (if a given knot is locally-flat, then $X$ is trivial). When considering both groups together, we sometimes refer to the ``knot group pair''. The following theorem generalizes Kervaire's \cite{Ke} necessary conditions for a group $\bG$ to be a knot group, $\pi_1(C)$.

\begin{theorem}\label{T: nec con}
The following conditions are necessary for the group $\bG$ to be the fundamental group of the complement, $C$, of a (not necessarily locally-flat) PL knot $K\subset S^n$, $n\geq 3$:
\begin{enumerate}
\item \label{J: finite}\label{I: finite} $\bG$ is finitely presentable,
\item \label{J: Z}\label{I: Z} $\bG/[\bG,\bG]\cong \Z$,
\item  \label{J: H2}\label{I: H2} $H_2(\bG)=0$,
\item \label{J: closure}\label{I: closure} There exist an element  $\bg\in \bG$ such that $\bG$ is the normal closure of $\bg$ (i.e. $\bg$ is of \emph{weight one}).
\end{enumerate}
\end{theorem}

We will refer to conditions \eqref{J: finite}-\eqref{J: closure} as the \emph{Kervaire conditions} \eqref{J: finite}-\eqref{J: closure} on a group.

\begin{proof}
The proof is a slight generalization of that of Kervaire \cite{Ke} for smooth knots:

As in \cite{Ke}, condition \eqref{J: finite} holds because $C$ is homotopy equivalent to the complement of the open regular neighborhood of the knot, and this is a finite simplicial complex; condition \eqref{J: Z} is due to $C$  being a homology circle by Alexander duality; and condition \eqref{J: H2} follows from $C$ being a homology circle and the Hopf exact sequence (see \cite{L77})
\begin{equation*}
\begin{CD}
\pi_2(C)@>\rho >>H_2(C)@>>>H_2(\pi_1(C))@>>>0,
\end{CD}
\end{equation*}
where $\rho$ is the Hurewicz homomorphism.

Condition \eqref{J: closure} requires the most modification, but again the basic idea is Kervaire's. We show that the adjunction of one relation to the group $\pi_1(C)$ will kill it; $\bg$ can then be taken as the relator. Equivalently, we show that attaching a disk to $C$ will create a simply-connected space. In particular, choose a point $c_0$ at which the knot is locally-flat. Then, locally, the regular neighborhood of $c_0$ in $S^n$ is isomorphic to  a 2-disk bundle in $S^n$ over a neighborhood of $c_0$ in $K$ and whose boundary circle bundle lies in $\bd C$. Let $Q$ be the 2-disk fiber with center $c_0$.  We show that $C\cup Q$ is simply connected.

Suppose that $\alpha$ is a curve representing an element of $\pi_1(C \cup Q)$. By PL approximation, we may assume that $\alpha$ is PL. Since $S^n$ is simply-connected for $n\geq 2$, there exists a map $F: D^2 \to S^n$ such that $F|_{S^1}=\alpha$. We may also assume $F$ to be a PL map into $S^n$, and by general position, we can assume (by applying a homotopy if necessary) that $F(D^2)$ intersects the knot only at a finite number of locally-flat points, $\{F(b_i)\}$, $b_i\in \text{int}(D^2)$ (since the dimension of the singular set must be $\leq n-4$ as  seen in our discussion of stratified pseudomanifolds in Section \ref{S: cons}). We can further assume, by further modifying $F$ if necessary, that $F$  maps a small disk $D^2_i$ around each $b_i$ homeomorphically onto a disk representing the fiber over $F(b_i)$ in the $2$-disk bundle  that is a  regular neighborhood of $F(b_i)$ in $S^n$. Now choose paths, $w_i$, in the knot $K$ from each $F( b_i)$ to $c_0$. By general position, we may assume that these paths are disjoint from each other (except at $c_0$) and from the singular set of the knot. We can now homotop $F$, using these paths, so that each $D^2_i$ contains a disk $E^2_i\subset \text{int}(D^2_i)$ such that $F(E^2_i)=Q$ and $F(D^2_i-E^2_i)\subset C$. Roughly speaking, since each path $w_i$ lies in the locally-flat part of the embedding, we can homotop $F(E_i)$  to $Q$ in a neighborhood of the path while keeping its boundary disjoint from the knot and then stretch $F(D^2_i-E^2_i)$ into the trace of the resulting homotopy on the boundary of $E^2_i$. Then, since this trace is disjoint from the knot, it can be pushed back into the complement of the regular neighborhood. Once this modification has been accomplished for all $i$, we see that $\alpha$ is in fact nullhomotopic in $C\cup Q$. Therefore, $C\cup Q$ is simply connected.
\end{proof}

For the case of knots with point singularities, we can generalize Kervaire's \cite{Ke} classification in higher dimensions to obtain a full classification of knot group pairs in dimensions $n\geq 6$. By the geometric arguments of Sections \ref{S: cons}, this is equivalent to classifying the groups of locally-flat disk knots together with those of their locally-flat boundary sphere knots.

\begin{theorem}\label{T: class}
Suppose $n\geq 6$. For the groups $\bG$ and $G$ to be the fundamental groups of the respective complements, $C$ and $X$,  of a locally-flat PL disk knot $J: D^{n-2}\into D^n$ and its boundary locally-flat sphere knot and for $\phi: G\to \bar G$ to be the homomorphism induced by inclusion, it is necessary and sufficient that $G$ and $\bar G$ satisfy the Kervaire conditions with elements of weight one $g\in G$ and $\bar g\in \bar G$ such that $\phi(g)=\bar g$. 
\begin{comment}
\begin{enumerate}
\item \label{I: finite} $G$ and $\bG$ are finitely presentable,
\item \label{I: Z} $G/[G,G]\cong \bG/[\bG,\bG]\cong \Z$,
\item  \label{I: H2} $H_2(G)=H_2(\bG)=0$,
\item \label{I: closure} There exist elements  $g\in G$ and $\bar g \in \bG$ and a homomorphism $\phi: G\to \bG$ such that 
\begin{enumerate}
\item $G$ is the normal closure of $g$,
\item $\bG$ is the normal closure of $\bar g$,
\item $\phi(g)=\bg$.
\end{enumerate}
\end{enumerate}
\end{comment}
(Note: our construction will in fact yield smooth knots, giving a slightly stronger realization theorem.)
\end{theorem}

\begin{proof}[Proof of necessity:] 

The necessity of the statements involving $G$ alone  follow from Kervaire's classification of higher dimensional locally-flat knot groups \cite{Ke}. Those involving $\bG$ alone follow from Theorem \ref{T: nec con} and the fact that $C$ is homotopy equivalent to the complement of the  non-locally-flat sphere knot given by adding the cone pair on the boundary to our knot pair $(D^n, J)$.

For the map condition, $\phi$ can be taken as the map on $\pi_1$ induced by inclusion. If we fix a \emph{simple meridian} of the boundary sphere knot in $X$, i.e. an embedded circle that bounds an embedded disk which intersects the knot only in a single point, the inclusion takes this meridian to a simple meridian of the disk knot. By Kervaire's theorem for sphere knots and the proof of Theorem \ref{T: nec con} above, the elements $g\in\pi_1(X)\cong G$ and $\bar g\in \pi_1(C)\cong \bG$ whose normal closures generate the groups can be represented by such meridians. Hence, by an appropriate choice of meridians (and basepoints), $\phi(g)=\bg$. 

\end{proof}

\begin{proof}[Proof of sufficiency]

We generalize the construction of Kervaire \cite{Ke} (see also Levine \cite[\S 9]{L77}).

We begin by constructing a CW complex that will serve as a blueprint for a handlebody construction. First, we construct two separate complexes, $P_2$ and $Q_2$, such that $\pi_1(P_2)=G$ and $\pi_1(Q_2)=\bG$. In fact, since $G$ and $\bG$ are finitely presented, we can take each complex to be the one point union of a set of circles representing generators together with a set of $2$-disks attached to represent the relations. Now, let $P_1$ and $Q_1$ be the $1$-skeleta of $P_2$ and $Q_2$, consisting of the one point unions of the circles representing generators $\alpha_i$ and $\beta_i$ of $G$ and $\bG$.  Each $\phi(\alpha_i)$ can be represented by some product of generators of $\bG$, and we use this to define a base-point preserving map from $P_1$ to $Q_1$. In other words, define the map on the circle representing $\alpha_i$ to be a representation of $\phi(\alpha_i)$ in $\pi_1(Q_2)=\bG$, which we can assume to lie in $Q_1$. Let $I_{\phi}$ denote the mapping cylinder of the induced map $P_1\to Q_2$, and let $T_2=I_{\phi}\cup P_2$, the quotient along the inclusion of $P_1$ into both $P_2$ and $I_{\phi}$. If we abuse notation and let $\phi$ also stand for the map $P_1\to Q_1$, then $T_2\sim_{h.e.} P_2\cup_{\phi} Q_2$. 

Notice that $\pi_1(T_2)\cong \bG$. In fact, $P_2\cup_{\phi} Q_2$ has only the circles representing the generators $\beta_i$ as $1$-cells, and the only $2$-cells are the disks representing the relations in $\bG$ and the images of the $2$-cells from $P_2$. But the boundary of each $2$-cell, say $D$, in $P_2$ represents the $0$ element of $G=\pi_1(P_2)$, and so under the map of $1$-cells induced by $\phi$, $\bd D$ must be mapped to a product of generators of $\bG$ which already bounds in $\bG$, because $\phi$ is a homomorphism. Therefore, in $P_2\cup_{\phi} Q_2$, the $2$-cells from $P_2$ introduce no new relations among the $\beta_i$, and $\pi_1(P_2\cup_{\phi} Q_2)\cong \pi_1(Q_2)\cong \bG$. 

Notice also that $T_2$ is a CW complex of dimension $2$. It can be obtained from the disjoint union of $P_2$ and $Q_2$ by adding a $1$-handle, $\gamma$, to connect the $0$-skeleta and then attaching $2$-cells whose attaching maps represent $\alpha_i\phi(\alpha_i)^{-1}$.

Next, we are going to need to modify the CW pair $(T_2,P_2)$ to a pair $(T,P)$ so that each of $T$ and $P$ are homology circles. This will be needed below.

The modules $H_1(T_2)$ and $H_1( P_2)$ are already as desired because $H_1(P_2)\cong G/[G,G]\cong \Z$ and $H_1( T_2)\cong \bG/[\bG, \bG]\cong \Z$.  Also $H_i( P_2)=H_i( T_2)=0$ for $i>2$ since $P_2$ and $T_2$ are $2$-dimensional complexes. Now consider $H_2( P_2)$ and $H_2(T_2)$. Since there are no $3$-cells in $T_2$ or  $P_2$, these are each free abelian groups, as they are the kernels of the boundary maps on the free abelian chain groups $C_2( T_2)$ and $C_2( P_2)$. Furthermore, the same is true of $H_2(T_2,P_2)$ as the kernel of the boundary map of the chain group $C_2(T_2, P_2)$, which is free abelian because $C_2( P_2)$ is a free direct summand of $C_2(T_2)$. Again since there are no $3$-cells, we have the exact sequence 
\begin{equation*}
\begin{CD}
0@>>> H_2(P_2) @>>> H_2( T_2) @>>> H_2(T_2, P_2) @>\bd_*>>,
\end{CD}
\end{equation*} 
which we can truncate as 
\begin{equation*}
\begin{CD}
0@>>> H_2( P_2) @>>> H_2( T_2) @>>> \text{ker}(\bd_*) @>>> 0.
\end{CD}
\end{equation*}
Since $ \text{ker}(\bd_*)$ is also free as a subgroup of a free group, this sequence splits and $ H_2( T_2)\cong H_2( P_2) \oplus K$, where $K\cong \text{ker}(\bd_*)$. Notice also that these free groups are all finitely generated since the chain groups are all generated  by cells of $T_2$, which is a finite complex by construction.

Let $F$ and $\bar F$ be free groups with the same ranks as $ H_2( P_2)$ and  $K$, respectively. Let $\{f_i\}$ and $\{\bar f_i\}$ represent generators of $F$ and $\bar F$, and let $\{e_i\}$ and $\{\bar e_i\}$ represent generators of $ H_2( P_2)$ and  $K$.  Let $\bd: F \oplus \bar F\to H_2( P_2) \oplus K $ be the homomorphism which takes each $f_i$ to $e_i$ and each $\bar f_i \to \bar e_i$. This map is clearly  injective, and furthermore the restriction $\bd|_F: F\to  H_2( P_2)$ is injective. We also claim that each element $e_i\in H_2(P_2)$ can be represented by a $2$-sphere in $P_2$ and that each element $\bar e_i\in H_2(T_2)$ can be represented by a $2$-sphere in $T_2$. These claims follow from the Hopf exact sequences
\begin{equation*}
\begin{CD}
\pi_2(P_2)@>\rho >>H_2(P_2)@>>>H_2(G)@>>>0\\
\pi_2(T_2)@>\rho >>H_2(T_2)@>>>H_2(\bG)@>>>0,
\end{CD}
\end{equation*}
since $ H_2(G)= H_2(\bG)=0$ by assumption. Therefore, let $\{E_i\}$ and $\{\bar E_i\}$ be $2$-spheres representing $\{e_i\}$ and $\{\bar e_i\}$ in $P_2$ and $T_2$. We attach $3$-cells to $T_2$ along the $E_i$ and $\bar E_i$. Let $T$ be the resulting complex, and let $P$ be the subcomplex complex obtained from attaching the $3$-cells along the $E_i$ to $P_2$.  Then $H_i( P)=H_i( T)=0$, $i\geq 2$, $H_1( T)\cong H_1( T_2)$, $H_1( P)\cong H_1( P_2)$.  Therefore, we can conclude that the CW $3$-complexes $P$ and $T$ are both homology circles.

We now indicate how to create a handlebody based upon the CW-pair $(T,P)$; we outline the procedure in this paragraph and provide the technical justifications in the next two. First, we construct a handlebody of dimension $n\geq 6$ modeled on $P$. In other words, begin with a ball $B^{n}$ and then attach one $1-$handle for each for each generator $\alpha_i$ and then $2$ and $3$ handles as prescribed by the construction of $P$.  Let $M$ denote this handlebody, which will be homotopy equivalent to $P$. Now consider $M\times I$, and add the remaining handles as handles of dimension $n+1$ to create  an $n+1$-manifold homotopy equivalent to $T$. We can assume that no  handles are attached to $M\times 0$.  
\begin{comment}\end{comment}
Let us call this new handlebody $N$; it will be homotopy equivalent to $T$.

We now justify this procedure by using some basic constructions from surgery theory. Let us first construct $M$ inductively. Suppose that we have chosen an order on the cells of $P$ so that $P$ is constructed by attaching the cells in order and with increasing dimensions. Suppose that $P_k$ is the subcomplex consisting of the first $k$ cells of the construction, and suppose that $D_{k+1}$ is the next cell to attach and that $D_{k+1}$ has dimension $r+1$. Suppose that we have an $n$-dimensional handlebody $M_k$ which consists of handles of core dimension $\leq 3$, that $M_k$ is stably parallelizable, and that there is a homotopy equivalence $f:M_k\to P_k$. We want to add a handle to $M_k$ to create $M_{k+1}$ which has the same properties and is homotopy equivalent to $P_{k+1}=P_k\cup D_{k+1}$. If $\alpha:\bd D_{k+1}\to P_k$  is the attaching map, the homotopy equivalence of $M_k$ and $P_k$ gives a homotopy class in $\pi_r(M_k)\cong\pi_r(P_k) $. Since $r\leq 2$, $n>5$, and $M_k$ consists of handles of core dimension $\leq 3$, any element of $\pi_r(M_k)$ can be represented by a homotopy element  in $[\tilde \alpha]\in\pi_r(\bd M_k)$ by general position. We want to form $M_{k+1}$ by adding an $r+1$ handle along $\tilde \alpha$. Then clearly we could obtain a homotopy equivalence $M_{k+1}\sim_{h.e.}P_{k+1}$ by extending $f$. The issue of whether we can attach such a handle is the standard surgery problem on $\bd M$ given by the diagram
\begin{equation*}\begin{CD}
S^r &@>\text{inclusion}>>& D^{r+1}\\
@VV\tilde{\alpha}V&& @VV\text{attachment}V\\
\bd M_k& @>f>>  & P_{k+1}.
\end{CD}\end{equation*} 
But since $r<n/2$, this surgery problem has a smooth solution, and the union of $M_k$ with the trace of the surgery gives us a new handlebody $M_{k+1}$ which has handles of the appropriate dimension and to which the stable parallelizablility extends (see \cite[Chapter 1]{W}, where we have take $\nu$ equal to the trivial bundle over $P_k$). We then continue by induction. Note that as the base step, we can use the trivial map of an $n$-ball to the base point of $P$.

Next, we want to add handles to $M\times I$ to obtain a handlebody $N$ homotopy equivalent to $T$. This time  $N$ will be the trace of a sequence of surgeries on the interior of $M$.  So let us order the cells of $T- P$, and suppose that we have created a stably parallelizable manifold $N_k$ homotopy equivalent to $T_k$ by adding a sequence of handles of core dimension $\leq 3$ to $M\times I$ in such a way as to leave $M\times 0$ invariant. The next attaching map $\alpha: S^r\to T_k$ gives an element of $\pi_r(N_k)$, and by general position, we may assume that this element can be represented by an element $[\td \alpha]\in \pi_r(\bd N_k)$. Furthermore,  by general position, we can assume that $\td\alpha(S^r)\subset \bd N_k- (M\times 0)$. Again we have a surgery problem, and the hypotheses on the dimensions and the stable parallelizability of $N_k$, and hence of $\bd N_k$, allow us to attach the next handle to create $N_{k+1}$ with the desired properties. The base case for this induction can be taken as the homotopy equivalence $M\times I\to P$ induced by the equivalence $M\to P$. After attaching all of the handles, we obtain $N$.

Having created $N$, we are really interested in the manifold $\bd N-\text{int}(M\times 0)$, whose boundary is $\bd M\times 0$. Let us call this manifold/boundary pair $(Z,Y)$. The boundary $Y$ is, in fact, that space which is created by Kervaire and then surgered to create a sphere knot with group $G$. In our procedure, this will be the result on the boundary. From Levine's treatment \cite{L77} of Kervaire's construction, we know that $\pi_1(Y)\cong \pi_1(P)\cong G$ and $H_*(Y)\cong H_*(S^1\times S^{n-2})$. We next claim that $\pi_1(Z)\cong \pi_1(T)\cong \bG$, that the inclusion $Z\into Y$ induces $\phi: G \to \bG$ on $\pi_1$, and that $H_*(Z)\cong H_*(S^1\times D^{n-1})$.

As an $n+1$-dimensional handlebody consisting of handles of core dimension $\leq 3$, $\pi_i(N, \bd N)=0$ for $i<n-2$ by general position. Similarly, $ \pi_i(\bd N, \bd N-\text{int}(M\times 0))\cong \pi_i(\bd N,Z)=0$ for $i<n-3$. Therefore, $\pi_i(Z)\cong \pi_i(N)$ (induced by inclusion) for $i<n-4$. Since $n>5$, this implies that $\pi_1(Z)\cong \bG$. To verify that the map on $\pi_1$ is as desired, consider the commutative diagram of inclusions
\begin{equation}\label{D: incs}
\begin{CD}
Y& @>>> &Z\\
@VVV&&@VVV\\
M&@>>>& N.
\end{CD}
\end{equation}
The vertical inclusion maps induce isomorphism on $\pi_1$; by Kervaire for the left map and by the above for the right map. Since the map $\pi_1(M)\to \pi_1(N)$ is the same as $\phi:\pi_1(P) \to \pi_1(T)$, the top map is also $\phi:G\to \bG$.

We next show that $Y\cong \bd M$ is a homology $S^1\times S^{n-2}$, and $\bd N$ is a homology $S^1\times S^{n-1}$. This follows just as in the proof of \cite{L77}[Lemma 10.1]; we briefly recall the reasoning: Since $\bd M$ is an $n-1$-manifold and $\bd N$ is an $n$-manifold, it suffices to prove that the homologies agree with those claimed in dimensions less than or equal to $\frac{n-1}{2}$ and $\frac{n}{2}$, respectively. 
As noted in the last paragraph, $(N, \bd N)$ is $n-3$ connected, and similarly, since $M$ is an $n$-dimensional handlebody consisting of handles of core dimension $\leq 3$, $(M, \bd M)$ is $n-4$ connected. Thus $H_i(\bd M)\cong H_i(M)\cong H_i(P)$ and $H_j(\bd N)\cong H_j(N)\cong H_j(T)$ for $i\leq n-5$ and $j\leq n-4$. Since $M$ and $N$ are homology circles, being homotopy equivalent to $P$ and $T$, we obtain the desired homology groups in these dimensions. The only remaining cases of concern are those where $n=6$, $i=2$, and $j=3$, or when $n=7$ and $i=3$. For $n=6$, we have the exact sequences
\begin{equation*}
\begin{CD}
H_3(M,\bd M)@>>> H_2(\bd M) @>>> 0\\
H_4(N,\bd N) @>>> H_3 (\bd N) @>>> 0,
\end{CD}
\end{equation*}
since $M$ and $N$ are homology circles. But  $H_3(M,\bd M)\cong H^3(M)=0$ and $H_4(N,\bd N)\cong H^2(N)=0$, by Lefschetz duality and  because $M$ and $N$ are homology circles. Therefore, $ H_2(\bd M)= H_3 (\bd N)=0$. Similarly, when $n=7$, we have 
 \begin{equation*}
\begin{CD}
H_4(M,\bd M)@>>> H_3(\bd M) @>>> 0,
\end{CD}
\end{equation*}
and we use $H_4(M, \bd M)\cong H^3(M)=0$.

Consider now the Mayer-Vietoris sequence
\begin{equation*}
\begin{CD}
@>>> H_i(Y) @>>> H_i(M\times 0) \oplus H_i(Z)@>>> H_i(\bd N) @>>>.
\end{CD}
\end{equation*}
%As regular neighborhoods of dimensions $>5$ of $P$ and $T$, which are complexes %of dimensions $3$ and homology circles, %

\noindent We know that $H_i(Y)\cong H_i(S^1\times S^{n-1})$ and $H_i(\bd N)\cong H_i(S^1\times S^n)$. Furthermore, $H_i(M)\cong H_i(P)\cong H_i(S^1)$. Therefore, around dimension $1$, the reduced exact sequence is 

\begin{equation*}
\begin{CD}
% 0@>>> &H_1(Y)& @>>>& H_1(M\times 0) \oplus H_1(Z)&@>>>& H_1(\bd N)& @>>>&0\\
% @V\cong VV&@V\cong VV&&@V\cong VV&&@V\cong VV&&@V\cong VV\\
0@>>> &\Z& @>>>& \Z \oplus H_1(Z)&@>>>& \Z&  @>>>& 0&,\\
\end{CD}
\end{equation*}
where we have used that $Z$ is connected. This sequence must split, and thus $H_1(Z)\cong \Z$. Also from the sequence, all other $H_i(Z)$ must be $0$ below and at the middle dimension. Since $Y$ is a homology $S^1\times S^{n-1}$, $Z$ will be a homology $S^1\times D^n$ if we show that the map induced by inclusion $H_1(Y)\to H_1(Z)$ is an isomorphism, since this will imply that $H_*(Z,Y)\cong 0$ below and at the middle dimensions, and we can then employ Lefschetz duality. 

To show that the map induced by inclusion $H_1(Y)\to H_1(Z)$ is an isomorphism, once again consider the diagram of inclusions \eqref{D: incs}.
As we have observed, the vertical inclusion maps induce isomorphism on $\pi_1$, so the homology diagram of the same spaces commutes with vertical isomorphisms. Therefore, it suffices to show that the map induced by inclusion $H_1(M)\to H_1(N)$ is an isomorphism, and, by homotopy equivalence, this is equivalent to showing that the map induced by inclusions $H_1(P)\to H_1(T)$ is an isomorphism. 

Recall that $g$ and $\bg$ are the defined as the elements of $G$ and $\bG$ whose normal closures are the whole groups $G$ and $\bG$. But then the  images of these elements also generate $G/[G,G]\cong \Z$ and $\bG/[\bG,\bG]\cong \Z$ (see \cite[\S 8]{L77}). Since, by hypothesis, $\phi$ takes $g$ to $\bg$, we see by the construction of $T$ that the inclusion of a circle representing $g$ in $P$ represents $\bg$ in $T$. In particular, this shows that the inclusion $P\into T$ induces an isomorphism on $H_1$. 

So, to summarize, at this point we have a manifold pair $(Z,Y)$ such that $\pi_1(Y)\cong G$, $\pi_1(Z)\cong \bG$, inclusion induces the map $\phi:G \to\bG$ of fundamental groups, $Y$ is a homology $S^1\times S^{n-2}$ and $Z$ is a homology $S^1\times D^{n-1}$. Finally, we attach a $2$-handle to $Z$ along a circle representing the element $g$ of $\pi_1(Y)$ (again, note that the dimension conditions and stable parallelizability are sufficient to add a handle). We claim that the resulting manifold/boundary pair $(\Delta, \Sigma)$ is an $n$-disk and that the cocore, $J\cong D^{n-2}$, of the handle is the knot which gives us the desired properties. 

Clearly, $\Delta-J$ is homotopy equivalent to $Z$, and therefore possesses the correct fundamental group. Similarly, $\bd \Delta-\bd J$ is homotopy equivalent to $Y-g$, letting $g$ also stand for an embedded $S^1$ in $Y$ representing $g\in \pi_1(Y)$. Therefore, by general position, $\pi_1(Y)\cong \pi_1(\bd \Delta-\bd D^2)$. The inclusion of these complements induces the appropriate maps. It remains to see that $\Delta$ is a ball. That $\Sigma=\bd \Delta$ is a sphere follows from Kervaire's proof of the classification theorem for sphere knot groups (see \cite{L77}, since that construction is the restriction of ours to the boundary). Furthermore, $\Delta$ is simply connected because we attached the handle to $g\in \pi_1(Y)$ which we know also represents $\bg\in \pi(Z)$. Since $\bG$ is the normal closure of $\bg$, this kills the fundamental group. Furthermore, since $Z$ was a homology circle and we have added a handle that kills a free generator of the homology in dimensions $1$,  $\td H(\Delta)=0$. It follows from the Whitehead theorem that $\Delta$ is contractible and bounded by a sphere. Since $n>5$, it follows from the higher dimensional Poincare conjecture that $\Delta$ is a PL ball, and the cocore of the last handle gives the desired knot. More generally, since handlebodies can always be assembled smoothly, we can conclude that we have a smooth knotted homotopy disk pair. But since a smooth homotopy sphere which bounds a homotopy ball is standard by the h-cobordism theorem, we actually obtain a smooth disk knot. The embedding respects the structure since it is the standard embedding within the last handle. 

\end{proof}

\begin{remark}
If  we add the hypotheses to the theorem that $G$ has a presentation with one more generator than relation and that $H_2([G,G])=0$, then there exists a knotted smooth homotopy disk pair $\Delta^3\subset \Delta^5$ with the given groups. For in this case, according to Levine \cite{L77}, the complex $P$ we constructed can be chosen to be $2$-dimensional. Then the dimensions work out sufficiently to apply  the above arguments to create manifolds $M$ of dimension $5$ and $N$ of dimension $6$. Then the arguments proceed as in the proof of the theorem, and we obtain a smooth knotted homotopy  $5$-disk pair with boundary a knotted homotopy $4$-sphere pair. Since the $5$-dimensional h-cobordism theorem holds in the topological category, we also see that we can obtain a topological disk knot pair in these dimensions with the desired groups. 

\end{remark}

\begin{corollary}\label{C: decouple}
Suppose $G$ and $\bG$ are two groups which satisfy the Kervaire conditions.
\begin{comment}
\begin{enumerate}
\item $G$ and $\bG$ are finitely presentable,
\item  $G/[G,G]\cong \bG/[\bG,\bG]\cong \Z$,
\item  $H_2(G)=H_2(\bG)=0$,
\item  There exist elements  $g\in G$ and $\bar g \in \bG$ such that $G$ is the normal closure of $g$ and $\bG$ is the normal closure of $\bar g$.
\end{enumerate}
\end{comment}
Then for $n\geq 6$, there exists a smooth disk knot $J:D^{n-2}\subset D^n$ such that the groups $\bG$ and $G$ are the fundamental groups of the respective complements, $C$ and $X$,  of the disk knot  and its boundary sphere knot. If  $G$ has a presentation with one more generator than relation and $H_2([G,G])=0$, then there exists a smooth knotted homotopy disk pair $\Delta^3\subset \Delta^5$ with this fundamental group pair.
\end{corollary}
\begin{proof}
This follows from the theorem and the above remark provided there exists a homomorphism $\phi:G\to \bG$ that takes $g$ to $\bg$. But we can always construct such a map as follows: Let $\psi:G\to \Z$ be the abelianization map, which we know must take $g$ to a generator. Assume that $g\to 1$ (otherwise compose with the isomorphism $\Z \to \Z$ which takes $1\to -1$). Then let $\eta: \Z \to \bG$ be the homomorphism uniquely determined by $1\to \bg$. Now take $\phi=\eta\psi$.  
\end{proof}

The realization part of the above theorem is still not completely satisfactory if we wish to make the transition back to singular sphere knots. In fact, if we take $G=\Z$ in the theorem, then the construction yields a disk knot whose boundary sphere knot is actually the unknot. If we then attach the cone on the boundary, we obtain a sphere knot pair, but one that is locally-flat everywhere since the cone on the trivial sphere knot pair yields the trivial disk knot pair. To fix  things up, we now indicate how to modify the above construction so that even if we desire $G=\Z$, we can obtain a non-trivial boundary knot.

The key is simply to use a different choice for the complex $P$ in the construction. Let us construct $P$ by first building its universal cover $\td P$. As the one $1$-skeleton of $\td P$, we take the real line with vertices on the integers and edges connecting these vertices so that there is one edge for each vertex pair $[k, k+1]$. To form the $2$-skeleton, we add a $2$-sphere $S^2$ at each vertex (in other words, we attach one  $2$-cell at each vertex using the trivial attaching map which takes $\bd D^2$ to the vertex). So far, we have a $2$-complex $\td P_2$ which clearly admits a free $\Z$ action. Furthermore, it is obvious that, making an arbitrary choice of basepoint,  $\pi_1(\td P_2)=0$ and $\pi_2(\td P_2)\cong H_2(\td P_2)\cong \Z[\Z]=\Z[t, t^{-1}]$ as an abelian group and as a $\Z$ module through the $\Z$ action induced via the free $\Z$ action of translations on the space. Choose now a Laurent polynomial $p(t)\in \Z[t,t^{-1}]$ such that $p(1)=1$ but $p(t)$ is not identically $1$. This polynomial represents an element of $\pi_2(\td P_2)\cong H_2(\td P_2)$ and so can be represented by the image of a sphere $S^2$ in $\td P_2$. To form $\td P$, we attach one $D^3$ along this sphere and along each of its translates under the $\Z$ action. Therefore, we have created a space $\td P$ which still admits a free $\Z$ action, $\pi_1(\td P)$ is still trivial, and $H_2(\td P)\cong \pi_2(\td P)\cong \Z[t,t^{-1}]/(p(t))$ by the construction and the Hurewicz theorem. Note also that $H_3(\td P)=0$, as $C_3(\td P)$ and $C_2(\td P)$ are generated by the $3$- and $2$-cells we attached and the boundary map $C_3(\td P)=\Z[\Z]\to C_2(\td P)=\Z[\Z]$ corresponds to multiplication by $p(t)$, which is an injective homomorphism ($\Z[\Z]$ is an integral domain).  

Now, let $P$ be the quotient space of $\td P$ under the $\Z$ action. Clearly $P$ is a space with one cell in each of dimensions zero through three and no other cells. By covering space theory, $\pi_1(P)\cong Z$, generated by the $1$ cell. Also $\pi_2(P)\cong \pi_2(\td P)\cong  \Z[t,t^{-1}]/(p(t))$. Let us compute the homology of $P$. The projection $\td P\to P$ is an infinite cyclic cover since it has covering $\Z$ action generated by $t$ and so there is a long exact sequence (Milnor \cite{M68})
\begin{equation*}
\begin{CD}
@>>> H_i(\td P) @>t-1>> H_i(\td P) @>>> H_i(P) @>>> H_{i-1}(\td P) @>>>,
\end{CD}
\end{equation*}
where the map $t-1: H_i(\td P)\to H_i(\td P)$ is that induced by the multiplication treating $H_i(\td P)$ as a $\Z[t, t^{-1}]$-module. We already know that $H_i(\td P)=0$ for $i\geq 3$ and $i=1$. This implies that $H_1(P)\cong \Z$ since clearly $H_0(\td P)\cong H_0(P)\cong \Z$ and $t$ acts trivially on $H_0(\td P)$. Furthermore, multiplication by $t-1$ is an automorphism of $ \Z[t,t^{-1}]/(p(t))$. Rather than construct an algebraic proof of this fact, we note that this result is well known in knot theory because we can construct knots whose complements $S^n-K$ are homology circles but whose Alexander modules in, say, dimension two (i.e. $H_2$ of the infinite cyclic cover of $S^2-K$) are $ \Z[t,t^{-1}]/(p(t))$ (see, e.g., Levine \cite{L66}). Thus $H_2(P)$ and $H_3(P)$ are $0$, so $P$ is a homology circle with non-trivial $\pi_2$ and infinite cyclic $\pi_1$.

The rest of the construction of a disk knot now goes through just as in the proof of the theorem using this $3$-complex $P$ instead of the obvious one. It remains to see that the boundary knot we obtain is non-trivial. Recall that the knot complement of the boundary knot in the construction is $\bd M$, where $M$ is a manifold of dimension $\geq 6$ that is homotopy equivalent to $P$. By general position with respect to the cores of the handles, the inclusion $\pi_i(\bd M)\to \pi_i(M)$ is surjective for $i< n-3$. Since $n\geq 6$, $\pi_2(\bd M)\to \pi_2(M)$ is surjective. But $\pi_2(M)$ is nontrivial, thus so is $\pi_2(\bd M)$. Therefore, $\bd M$ can not be the complement of the trivial knot as the complement of the trivial knot is  homotopy equivalent to $S^1$.

As a side note, observe that if $n\geq 7$ then $\pi_2(\bd M)\cong \pi_2 (M)$ and $\pi_2(\bd M)\cong \pi_2( \widetilde{\bd M})\cong H_2(\widetilde{\bd M})$. Hence in this case the boundary knot will have $p(t)$ as its Alexander polynomial in dimension two.

We summarize this discussion as a proposition:

\begin{proposition}\label{P: nontriv}
Given a pair of groups satisfying the hypotheses of Theorem \ref{T: class} or Corollary \ref{C: decouple}, there exist locally-flat disk knots which satisfy the conclusions of the theorem or corollary and whose boundary sphere knots are not the trivial knot (and hence these disk knots are also not trivial). 
\end{proposition}

\section{Boundary knot groups}\label{S: links}

For sphere knots with singular sets of dimension greater than $0$, the necessary conditions of Theorem \ref{T: nec con} might no longer hold for the boundary knot. In particular, calculations below will show that $H_1(X)\cong H_1(\Sigma)\oplus \Z$, so the boundary knot group may not abelianize to $\Z$.
Condition \eqref{I: closure} may also fail to hold. For example, if we frame spin a knot $K$ with a point singularity around a manifold $M$ (see \cite{GBF}, \cite{GBF1}, and below) and $G$ is the boundary knot group of $K$, then the boundary knot group of the spun knot $\sigma(K)$ will be $G\times \pi_1(M)$. Therefore, we can not, in general, expect the closure condition to hold since the elements of $G$ and $\pi_1(M)$ commute in $G\times \pi_1(M)$ and $\pi_1(M)$ may not have weight one. 

In order to be able to say something about boundary knot groups, we begin by calculating the homology of $X$, a calculation of independent interest which may prove useful in calculating other invariants of knots.

\begin{theorem}\label{T: H of X} For a PL knot $K\subset S^n$, 
\begin{equation*}
H_i(X)\cong
\begin{cases}
H_0(\Sigma), & i=0,\\
H_{i-1}(\Sigma)\oplus  H_{i}(\Sigma),& 1\leq i\leq n-3,\\
0, & i>n-3.
\end{cases}
\end{equation*}
(Note that $H_{n-3}(\Sigma)=0$ since $\Sigma$ can have dimension at most $n-4$.)
\end{theorem}
\begin{proof}
This formula holds trivially for $n\leq 3$, since in this case there can be no singular set of the PL embedding. For $n=4$, the singular set will be a set of isolated points, and $X$ will be a disjoint union of homology circles, one for each isolated singular point. Therefore, we can concentrate on the cases $n\geq 5$.

We will use the long exact sequence of the pair $(\bd C,X)$, so first we calculated the homologies of $\bd C$ and of the pair.
We first calculate the homology of $\bd C$. By Alexander duality, $C$ is a homology circle.  So, applying 
Poincare-Lefschetz duality, the universal coefficient theorem, and the long exact sequence of the pair, we obtain that 
\begin{equation*}
H_{i}(\bd C)\cong
\begin{cases}
\Z, &i=0,1, n-1,n-2,\\
0,& \text{otherwise.}
\end{cases}
\end{equation*}

\begin{comment}
We first calculate the homology of $\bd C$. By Alexander duality, $C$ is a homology circle, so
\begin{equation*}
H_i(C)\cong
\begin{cases}
\Z, &i=0,1,\\
0,& \text{otherwise,}
\end{cases}
\end{equation*}
and by Poincare-Lefschetz duality,
\begin{equation*}
H^{i}(C, \bd C)\cong
\begin{cases}
\Z, &i=n,n-1,\\
0,& \text{otherwise.}
\end{cases}
\end{equation*}
Therefore, we see from the universal coefficient theorem and the long exact sequence of the pair that 

\begin{equation*}
H_{i}(\bd C)\cong
\begin{cases}
\Z, &i=0,1, n-1,n-2,\\
0,& \text{otherwise.}
\end{cases}
\end{equation*}
\end{comment}

For the homology of the pair $(\bd C, X)$, let us denote by $T$ the manifold  $\bd C-\text{int}(X)$ with $\bd T\cong \bd X$. As the appropriate part of  the boundary of the regular neighborhood of the locally-flat part of the knot, $T$ is a tube homotopy equivalent to $S^1\times (K-\Sigma)$. The triviality of the bundle is ensured by the existence of the map $S^n-K\to S^1$ that takes any meridian of the knot to the circle with degree $1$. This map exists because $S^n-K$ is a homology circle; see \cite[Lemma 4.1]{GBF} or \cite{GBF1} for more details. By excision, Poincare-Lefschetz duality, and homotopy equivalence, $H_i(\bd C, X)\cong H_i(T,\bd T)\cong H^{n-1-i}(T)\cong H^{n-i-1}(S^1\times (K-\Sigma))$. This last group is $H^{n-i-1}(K-\Sigma)\oplus H^{n-i-2}( K-\Sigma)$ by the K\"unneth theorem. Finally, noting that $K-\Sigma$ is connected and applying Alexander duality, we obtain

\begin{equation*}
H_i(\bd C, X)\cong
\begin{cases}
\td H_{i-2}(\Sigma)\oplus \td H_{i-1}(\Sigma), &i<n-2,\\
H_{n-4}(\Sigma)\oplus \Z,&i=n-2,\\
\Z,&i=n-1,\\
0,&i>n-1.
\end{cases}
\end{equation*}

Now, since $H_i(\bd C)=0$ in the middle dimensions, the long exact sequence of the pair $(\bd C, X)$ shows us that $H_i(X)\cong H_{i+1}(\bd C, X)\cong \td H_{i-1}(\Sigma)\oplus \td H_{i}(\Sigma)$ for $2\leq i\leq n-4$. It remains to examine the ends of the sequence. The claim at bottom is easy to check, using the correspondence between connected components of $X$ and $\Sigma$. 

\begin{comment}
At the  bottom, we have
{\footnotesize \begin{equation*}
\begin{CD}
H_2(\bd C)&@>>> &H_2(\bd C, X)&@>>> &H_1(X)&@>>> &H_1(\bd C)&@>>> &H_1(\bd C, X)&@>>> &\td H_0(X)&@>>> &0\\
@VV\cong V&&@VV\cong V&&@VV\cong V&&@VV\cong V&&@VV\cong V&&@VV\cong V&&@VV\cong V&\\
0&@>>> &\td H_{0}(\Sigma)\oplus \td H_{1}(\Sigma)   &@>>> &H_1(X)&@>>> &\Z &@>>> &\td H_0(\Sigma)&@>d>> &\td H_0(X)&@>>> &0.\\
\end{CD}
\end{equation*}} 

\noindent Since components of $X$ correspond to components of $\Sigma$, the map $d$ must be an isomorphism, and so $H_1(X)\cong \td H_{0}(\Sigma)\oplus \td H_{1}(\Sigma) \oplus \Z\cong H_0(\Sigma)\oplus H_1(\Sigma)$.
\end{comment}

At the top of the sequence, since $X$ is an $n-1$ manifold with boundary, we have
\begin{comment}
{\footnotesize\begin{equation*}
\begin{CD}
H_{n-1}(X)&@>>>& H_{n-1}(\bd C)&@>d >> &H_{n-1}(\bd C, X)&@>>> &H_{n-2}(X)&@>>> &H_{n-2}(\bd C)&@>e>>\\
@VVV&& @VV\cong V&&@VV\cong V&&@VV\cong V&&@VV\cong V\\
0&@>>>&\Z &@>d>> &\Z &@>>> &H_{n-2}(X)&@>>> &\Z &@>e>>\\
&&&&&@>e>> &H_{n-2}(\bd C, X)&@>>> &\td H_{n-3}(X)&@>>> &H_{n-3}(\bd C)\\
&&&&&&&&@VV\cong V&&@VV\cong V&&@VV\cong V&\\
&&&&&@>e>> &\td H_{n-4}(\Sigma)\oplus\Z &@>>> &\td H_{n-3}(X)&@>>> &0,
\end{CD}
\end{equation*} }
\end{comment}

{\footnotesize\begin{equation*}
\begin{CD}
H_{n-1}(X)&@>>>& H_{n-1}(\bd C)&@>d >> &H_{n-1}(\bd C, X)&@>>> &H_{n-2}(X)&@>>> &H_{n-2}(\bd C)\\
% @VVV&& @VV\cong V&&@VV\cong V&&@VV\cong V&&@VV\cong V\\
% 0&@>>>&\Z &@>d>> &\Z &@>>> &H_{n-2}(X)&@>>> &\Z &@>e>>\\
&@>>> &H_{n-2}(\bd C)&@>e>> &H_{n-2}(\bd C, X)&@>>> &\td H_{n-3}(X)&@>>> &H_{n-3}(\bd C)\\
% &&&&&&&&@VV\cong V&&@VV\cong V&&@VV\cong V&\\
%&&&&&@>e>> &\td H_{n-4}(\Sigma)\oplus\Z &@>>> &\td H_{n-3}(X)&@>>> &0,
\end{CD}
\end{equation*} }

\noindent in which we know that the first and last groups are $0$, the second, third, and fifth are isomorphic to $\Z$, and $H_{n-2}(\bd C, X)\cong \td H_{n-4}(\Sigma)\oplus\Z$. We claim that $d$ and $e$ are both split injections, and this will suffice to finish the proof. To this end, consider the map $$H_1(T)\cong H_1(S^1)\oplus H_1(K-\Sigma)\to H_1(\bd C)\cong \Z$$induced by inclusion and in which we have utilized the isomorphisms $H_1(S^1)\otimes H_0 (K-\Sigma)\cong H_1(S^1)\otimes \Z\cong H_1(S^1)$ and similarly for the other term. From the K\"unneth theorem, the generator of the first summand of  $H_1(T)$ is a meridian of the knot around a locally-flat point of the embedding, and similarly for $\bd C$ due to the isomorphism $H_1(\bd C)\cong H_1(C)$ induced by inclusion. Hence this map is a surjection and a split surjection, since $\Z$ is free. By excision and the naturality of Poincare-Lefschetz duality, 
\begin{comment}
we obtain a diagram 
\begin{diagram}
&&H_1(T)&\rTo &H_1(\bd C)\\
&&\dTo^{\cong} &&\dTo^{\cong}\\
H^{n-2}(\bd C, X)&\rTo^{\cong}_{\text{excision}} & H^{n-2}(T,\bd T)&\rTo &H^{n-2}(\bd C)
\end{diagram}
(where, technically, we have used the homotopy equivalence of $T$ and $T-\bd T$), and so 
\end{comment}
we obtain a split surjection $H^{n-2}(\bd C, X)\cong  H^{n-2}(T,\bd T)\to H^{n-2}(\bd C)$ induced by the standard inclusion $C^*(\bd C, X)\into C^*(\bd C)$. Dually, this induces a split injection $\Hom(H^{n-2}(\bd C),\Z)\to \Hom( H^{n-2}(\bd C, X), Z)$. Notice now that $H^{n-1}(\bd C)\cong H_0(\bd C)\cong \Z$  is free, and so is $H^{n-1}(\bd C,X)\cong H_0(T)\cong \Z$. Therefore, the natural commutative diagram of universal coefficient sequences for $H_{n-2}(\bd C)$ and $H_{n-2}(\bd C, X)$ reduces to
\begin{comment}
\begin{equation*}
\begin{CD}
0&@<<<&\Hom(H^{n-2}(\bd C) , \Z) &@<<<&  H_{n-2}(\bd C)   &@<<<&\ext(H^{n-1}(\bd C) , \Z)  &@<<<&0\\
@VVV&&@VVV&&@VVV&&@VVV&&@VVV&&\\
0&@<<<&\Hom(H^{n-2}(\bd C,X) , \Z) &@<<<&  H_{n-2}(\bd C,X)   &@<<<&\ext(H^{n-1}(\bd C,X) , \Z)  &@<<<&0
\end{CD}
\end{equation*}
reduces to 
\end{comment}

\begin{equation*}
\begin{CD}
\Hom(H^{n-2}(\bd C) , \Z) &@<\cong<<&  H_{n-2}(\bd C)  \\
@VVV&&@VVV\\
\Hom(H^{n-2}(\bd C,X) , \Z) &@<\cong<<&  H_{n-2}(\bd C,X) . 
\end{CD}
\end{equation*}
Since the lefthand map is a split injection, so must be the righthand map, but this is exactly the standard map induced by projection and so the same map that we considered in the long exact sequence of $(\bd C, X)$. 

The claim for the map $H_{n-1}(\bd C)\to H_{n-1}(\bd C, X)$ is even simpler to prove since this is the canonical map of orientation classes.

\begin{comment}
 Equivalently, the Poincare-Lefschetz duality diagram in this case gives us
\begin{equation*}
\begin{CD}
&&H_0(T)&@>>> &H_0(\bd C)\\
&&@V\cong VV&&@V\cong VV\\
H^{n-1}(\bd C, X)&\cong& H^{n-1}(T,\bd T) &@>>> &H^{n-1}(\bd C),
\end{CD}
\end{equation*}
and the top map is the obvious isomorphism. The rest of the above argument then applies again since $H^n(\bd C)=H^n(\bd C, X)=0$ by dimension considerations.
\end{comment}
\end{proof}

\begin{corollary}\label{C: link group}
If $\Sigma$ is connected and $G$ is the group of the boundary knot, i.e. $G\cong \pi_1(X)$, then $G/[G,G]\cong H_1(X)\cong H_1(\Sigma)\oplus H_0(\Sigma)$, and there is a surjection $H_2(X)\cong H_2(\Sigma)\oplus H_1(\Sigma)\to H_2(G)$. 
\end{corollary}
\begin{proof}
The first claim follows from the theorem by abelianizing $G$ and the second from the Hopf exact sequence.
\end{proof}

\begin{corollary}
If $\Sigma$ is $2$-connected, then $G\cong \pi_1(X)$ and $\bG\cong \pi_1(C)$ must satisfy conditions \eqref{I: finite}-\eqref{I: closure} of Theorem \ref{T: class}.
\end{corollary}
\begin{proof}
The necessity of the conditions involving $\bG$  are the content of Theorem \ref{T: nec con}. Conditions \eqref{I: finite} is true for $G$ because $X$ has a finite complex as a deformation retract. Conditions \eqref{I: Z} and \eqref{I: H2} follow for $G$ from the previous corollary. The proof of condition \eqref{I: closure} for $G$ follows just as in the proof of the equivalent condition in Theorem \ref{T: nec con} once we observe that the pair $(\bd \bar{N}(\Sigma), \bd \bar{N}(\Sigma)\cap K)$ is a locally-flat pair and that $\pi_1(\bd \bar{N}(\Sigma))\cong \pi_1(\Sigma)$ by general position, since $\text{dim}(\Sigma)\leq n-4$. 

Finally, the map of condition \eqref{I: closure} is that induced by inclusion, which takes a simple meridian of the boundary knot into a simple meridian of the big knot. But we know that we can take these meridians to correspond to the elements $g$ and $\bar g$ whose normal closures give the whole groups.   
\end{proof}

Using Theorem \ref{T: class} and this corollary, we can state the following theorem:

\begin{theorem}\label{T: bd kerv}
Conditions \eqref{I: finite}-\eqref{I: closure} of Theorem \ref{T: class} on the group pairs $(\bG, G)$ are necessary and sufficient for these groups to be knot group pairs in the class of PL-knots which have $2$-connected singular sets. 
\end{theorem}

Of course this statement should not be read to imply that we can necessarily realize a group pair once we have fixed the singular set $\Sigma$; we only know how to construct such a knot with point singularities realizing the group pair. 

As noted above, if we drop these hypotheses on $\Sigma$, then clearly the conditions of Theorem \ref{T: class} may not hold for $G$ (though we have seen above in Theorem \ref{T: nec con} that they will still hold for $\bG$). However, we can say a little bit more about the boundary knot group $G=\pi_1(X)$ in another special case: when the knot can be stratified (see Section \ref{S: cons}) so that the singular set of the knot $K$ has only one connected stratum which is thus a connected manifold. In particular, $\pi_1(X)$ will be part of an exact sequence involving the fundamental group of the manifold singular set $\Sigma$ and the knot group of the link pair of $\Sigma$. 

It follows from the definition of a stratified pseudomanifold (see Section \ref{S: cons}) that if  the singular set, $\Sigma$, of a manifold embedding has only one stratum, then $\Sigma$ must be a  manifold of some dimension $n-k-1$ with $k\geq 3$. If $\Sigma$ is furthermore connected, every point of $\Sigma$ will have the same \emph{link pair} $(L, \ell)$, a locally-flat pair of spaces such that $L$ is PL-homeomorphic to $S^k$ and $\ell$ is PL-homeomorphic to $S^{k-2}$.  The pair $(L,\ell)$ can thus be regarded as representing a knot,  the \emph{link knot}.  Similarly, if a knot can be stratified so that $\Sigma$ possesses a connected component that intersects only one stratum (of dimension $<n-2$), then that component will be a manifold and have associated to it a unique link knot.

\begin{theorem}\label{T: manif sing}
Let $(S^n, K)$ be a PL knot. Suppose that $(S^n, K)$ can be stratified so that some connected component of the singular set intersects only one stratum. Let $M^{n-k-1}$ denote this manifold  component and  $(L, \ell)$ its link knot. Let $X$ be the exterior of the boundary knot corresponding to this component. Then there is a long exact sequence 
\begin{equation*}
\begin{CD}
 @>>> \pi_i(L-\ell)@>>> \pi_i(X) @>>> \pi_i(M)@>>>\pi_{i-1}(L-\ell) @>>> .
\end{CD}
\end{equation*}
In particular, if $G=\pi_1(X)$ is the boundary knot group, then
\begin{equation*}
\begin{CD}
@>>> \pi_1(L-\ell)@>>> G @>>> \pi_1(M)@>>> 0
\end{CD}
\end{equation*}
is exact. Furthermore, if $\lambda$ is an element of the knot group $\pi_1(L-\ell)$ whose normal closure is the whole group, then the image of $\lambda$ in $G$ is an infinite cyclic subgroup.
\end{theorem}
\begin{proof}
Recall that, by definition, $X$ is the exterior of the knot $K$ in the boundary of a closed regular neighborhood $\bar N(M)$ of $M$. But this neighborhood is stratum-preserving homotopy equivalent to the mapping cylinder of a stratified fibration $p: E\to M$ (see \cite{GBF} or  \cite{GBF3} for proof and explanation of the terminology). Furthermore, the stratified fibration $p: E\to M$ has as stratified fiber a stratified space that is stratum-preserving homotopy equivalent to the pair $(L,\ell)$. In particular, this implies that the restriction of  $p$ to the top stratum of $E$, say $\mc E$, is a fibration whose fiber, $\mc{L}$, is homotopy equivalent to  $L-\ell$. Therefore,  from this fibration we obtain a long exact sequence 
\begin{equation*}
\begin{CD}
@>>> \pi_i(L-\ell)@>>> \pi_i(\mc E) @>>> \pi_i(M)@>>>.
\end{CD}
\end{equation*}

It remains to show that   $\pi_i(\mc E)\cong \pi_i(X)$. Let $Y$ denote the mapping cylinder of $p:E\to M$ (stratified so that $M$ is the bottom stratum with the higher dimensional skeleta being the mapping cylinders of those of $E$). Then $E$ is stratum-preserving homotopy equivalent to $Y-M$ and, in particular, $\mc{E}$ is homotopy equivalent to the top stratum of $Y-M$. But $Y$ is stratum preserving homotopy equivalent to $\bar N(M)$ so that $\mc{E}\sim_{h.e.} \bar N(M)-\bar N(M)\cap K$. Finally, we note that the inclusion of $X$ into $\bar N(M)-\bar N(M)\cap K$ is a homotopy equivalence. This follows since both spaces are triangulable and the inclusion is a weak homotopy equivalence by an easy corollary to  the stratified generalized annulus property for regular neighborhoods (see the proof of \cite[ Lemma 8.7]{GBF2}; this lemma actually concerns intersection homology, but the proof applies  equally well to the homotopy groups of a fixed stratum).

To prove the last statement of the theorem, we can consider $L$ as a subset of $\bar N(M)-M$ by choosing a sufficiently small distinguished neighborhood of some point of $M$ and choosing a homeomorphic image of $L$ determined by the product structure of the distinguished neighborhood. As a consequence of the discussion in the proof of \cite[Theorem 6.1]{GBF3}, the following diagram exists and commutes up to stratum-preserving homotopy:
\begin{equation*}
\begin{CD}
L  &@>>>&  \mc{L}\\
@VVV&&@VVV\\
\bar N(M) &@<<<& Y,
\end{CD}
\end{equation*}
where $\mc{L}$ is a fiber of $p:E\to M$, the vertical maps are inclusions, and the top and bottom maps are stratum-preserving homotopy equivalence. This implies, in particular, that the map $\pi_1(L-\ell)\to G$ of our long exact sequence is, up to isomorphism, that induced by an inclusion of $L-\ell$ into $\bar N(M)-K$ (taking into account the homotopy equivalence of $X$ and $\bar N(M)-K$). 
Now consider the following commutative diagram in which the vertical maps represent abelianizations and the horizontal maps are induced by inclusions: 
\begin{equation}\label{D: abel}
\begin{CD}
\pi_1( L-\ell) &@>>>& G\cong \pi_1(\bar N(M)-K)\\
@VVV&&@VVV\\
H_1(L-\ell)\cong Z& @>>> & H_1(\bar N(M)-K)&.
\end{CD}
\end{equation}
Up to isomorphism, the map of our exact sequence is that represented by the top line of this diagram. So let $\lambda$ be an element of $\pi_1(L-\ell)$ whose normal closure is the whole link knot group.  We know from our previous discussions that the left map takes $\lambda$ onto a generator. But any cycle representing the image of $\lambda$ in $H_1(L-\ell)$ is homologous to a  simple meridian of the knot  $\ell$ (up to orientation), and if any multiple of the image of such  a cycle  bounded in $\bar N(M)-K$ then it would also bound in $S^n-K$. But this is impossible since a simple meridian of $\ell\subset K$ is also clearly a simple meridian of $K$ which we know generates $H_1(S^n-K)\cong Z$. Thus the composite $H_1(L-\ell)\to H_1(\bar N(M)-K)\to H_1 (S^n-K)$ is an isomorphism and the image of $\lambda$ must generate an infinite cyclic subgroup of $H_1(\bar N(M)-K)\cong H_1(X)$. By the commutativity of the diagram, it follows that the image of $\lambda$ must also generate such a subgroup in $G$.
\end{proof}

\begin{comment}
%consider the following diagram in which both vertical arrows are abelianization %maps and the bottom map is induced by inclusion:
%\begin{equation}\label{D: abel}
%\begin{CD}
%\pi_1( L-\ell) &@>>>& G\cong \pi_1(X)\\
%@VVV&&@VVV\\
%H_1(L-\ell)cong Z& @>>> & H_1(X).
%\end{CD}
%\end{equation}
%We claim that this diagram commutes. In fact, 
we can consider $L$ as a subset of $\bar N(M)-M$ by choosing a sufficiently small distinguished neighborhood of some point of $M$ and choosing a homeomorphic image of $L$ determined by the product structure of the distinguished neighborhood. Then clearly we have commutativity of the diagram
\begin{equation}\label{D: abel}
\begin{CD}
\pi_1( L-\ell) &@>>>& \pi_1(\bar N(M)-M) &@<\cong <<&\pi_1(X)\\
@VVV&&@VVV&@VVV\\
H_1(L-\ell)\cong Z& @>>> & H_1(\bar N(M)-M)&@<\cong <<& H_1(X)
\end{CD}
\end{equation}
in which all horizontal maps are induced by inclusions and all vertical maps are abelianizations. To see that the induced map $\pi_1(L-\ell)\to \pi_1(X)$ is the same as that determined by the long exact homotopy sequence of the fibration it is sufficient to know that the homotopy equivalence $\bar N(M) \to Y$ actually restricts to a homotopy equivalence from $L$ to the fiber of $p$. This will be true is we make a convenient choice of the retracting function used to define the homotopy equivalence, in particular one which collapses $L$ along the cone lines in a distinguished neighborhood; see the proof of Theorem 6.1 in \cite{GBF3} for the details of these homotopy equivalences).
\end{comment}

\begin{corollary}\label{C:  s.c. man'}
In the setting of Theorem \ref{T: manif sing}, if $\Sigma$ is a simply-connected manifold, then there is a surjection $\pi_1(L-\ell)\to G$. If $\Sigma $ is $2$-connected, $\pi_1(L-\ell)\cong G$. 
\end{corollary}
\begin{proof}
This is clear from Theorem \ref{T: manif sing}.
\end{proof}

\begin{corollary}\label{C: s.c. man}
If $K$ is a PL knot possessing a stratification such that the singular set is a single stratum consisting of a $2$-connected manifold $\Sigma$, then the group of the boundary knot, $G=\pi_1(X)$, must satisfy the Kervaire conditions. 
\end{corollary}
\begin{proof}
Note that we already know that this corollary is true if $\Sigma$ is a point, so we can assume that $\Sigma$ has dimension $\geq3$ and the knot dimension $n\geq 7$ since dim$(\Sigma)\leq n-4$ and the hypotheses cannot hold if $\Sigma$ is a manifold of dimension $1$ or $2$.

From the previous corollary, we know that $G$ is the quotient of a knot group $\pi_1(L-\ell)$. It follows immediately that $G$ must be finitely generated since this is true of $\pi_1(L-\ell)$. Next, since it follows from Theorem \ref{T: H of X} that $H_1(X)\cong \Z$ and $H_2(X)=0$, it must be that $G$ abelianizes to $\Z$ and that $H_2(G)=0$ by the Hopf exact sequence. It remains to see that $G$ is the normal closure of one of its elements. This is equivalent to showing that there is an element $g\in G$ such that adding the relation $g=1$ to the presentation kills the group. But we already know that there exists such and element, say $x$, in $\pi_1(L-\ell)$ and that, since $G$ is a quotient of $\pi_1(L-\ell)$, $G$ can be presented by adding extra relations to a presentation of $\pi_1(L-\ell)$. Then the group presented by the presentation of $\pi_1(L-\ell)$ together with both the extra relations $x=1$ and those that must be added to obtain a presentation of $G$ is the trivial group. But clearly this is the same group obtained from the presentation of $G$ together with the relation given by setting the image of $x$ equal to $1$. In other words, $G$ is the normal closure of the image under the projection of any element of $\pi_1(L-\ell)$ whose normal closure is all of $\pi_1(L-\ell)$. We note for later use that, as usual, such an element can be represented by a simple meridian of the knot $\ell$ which is also homotopic to simple meridians of $K$ and of the boundary knot (assuming a fixed embedding of $(L, \ell)$ in $(S^n, K)$). 
\end{proof}

We will see below that for $n$ sufficiently large these conditions are also sufficient to classify boundary knot groups of knots with a single manifold stratum  if, in addition, the manifold can be embedded with a framing into a sphere. Hence we can obtain a complete characterization of the knot group pairs of knots of this type in a large range of dimensions.

\section{Realizing groups pairs for higher dimensional singular sets}\label{S: realize}

In order to demonstrate some actual examples of group pairs for knots with higher dimensional singular sets, let us calculate the knot groups of frame twist-spun knots.

First, we briefly review the definition of frame twist-spinning as given in \cite{GBF} and \cite{GBF1}. This construction generalizes the frame spinning of Roseman \cite{Ro89} and the twist-spinning of Zeeman \cite{Z65}. For a more detailed description, see \cite{GBF1} or \cite{GBF7}.

To set up the proper language, we adopt
some notation from Section 6 of Zeeman's paper, \cite{Z65}, in which he
introduces twist spinning. If we consider the unit sphere
$S^{m-1}$ in the Euclidean space $\R^m=\R^{m-2}\times \R^2$, then we can
define the latitude for a point $y\in S^{m-1}$ as its projection onto
$\R^{m-2}$ and its longitude as the angular polar coordinate of the
projection of $y$ onto the $\R^2$ term. Hence the latitude is always
well-defined, while the longitude is either undefined or a unique point of
$S^1$ dependent on whether or not $y$ lies in the sphere $S^{m-3}$ that is
the intersection of $S^{m-1}$ with $\R^{m-2}\times 0$. Notice that in the
case where the longitude in undefined, the point on the sphere is uniquely
determined by its latitude (just as on a standard globe). As in Zeeman's
paper, to simplify the notation in abstract cases, we will simply refer to the
latitude-longitude coordinates $(z, \theta)$ in either case.

Let $(D_{-}^{m}, D_{-}^{m-2})$ be an unknotted open disk pair which is
the open neighborhood pair of a point that does not lie in the singular set  of the
embedding of a knot $K\subset S^{m}$. Let $(D_{+}^{m}, K_+)=(S^{m}, K)-
(D_{-}^{m},
D_{-}^{m-2})$. This is a disk knot, possibly not locally-flat, with the unknotted
locally-flat sphere
pair as boundary. We can identify this trivial boundary
sphere pair $(S^{m-1}, S^{m-3})$ with the unit sphere in $\R^m$ and its intersection with $\R^{m-2}$. Using this identification, we can assign latitude-longitude coordinates $(z,\theta)$ to $S^{m-1}$ such that $S^{m-3}$ is the sphere with undefined longitude. If $M^k$ is a closed connected manifold then $M^k\times (D_{+}^{m},
K_+)$ gives a bundle of knots, and the points in $\bd[M^k\times
D_{+}^{m}]$ can be given coordinates  $(x,z,\theta)$, where $x\in M$ and
$(z,\theta)$ are the latitude-longitude coordinates of $\bd D^m$. 

Similarly, given an embedding of $M^k\subset S^{m+k-2}$ with framing $\psi$,
where $S^{m+k-2}$ is the $(m+k-2)$-sphere embedded in $S^{m+k}$ with the
standard normal bundle, we form
\begin{equation*}
(S^{m+k}, S^{m+k-2})-M^k\times
\text{int}(D^{m-2}\times D^2, D^{m-2}).
\end{equation*}
Again the boundary can be identified as $M^k \times (S^{m-1},
S^{m-3})$, and the framing $\psi$, together with the trivial framing of
$S^{m+k-2}$ in $S^{m+k}$, allows us to assign to this boundary the same
$(x,z,\theta)$-coordinates.

Given a map $\tau: M^k\to S^1$, we can form the frame twist-spun knot
$\sigma_M^{\psi,\tau}(K)$ of ambient dimension $n=m+k$ as
\[
[(S^{m+k}, S^{m+k-2})-M^k\times \text{int}(D^{m-2}\times D^2, D^{m-2})]
\cup_f [M^k\times (D_{+}^{m}, K_+)],
\]
where $f$ is the attaching homeomorphism of the boundaries 
\[f:\bd[M^k\times (D_{+}^{m}, K_+)]\to\bd[(S^{m+k}, S^{m+k-2})-M^k\times \text{int}(D^{m-2}\times D^2, D^{m-2})]\]
  which, identifying each with $M^k \times (S^{m-1},
S^{m-3})$ as above, takes
$(x,z,\theta)\to (x, z, \theta + \tau(x))$, where we define the addition in
the last coordinate as the usual addition on $S^1=\R/\Z$. The map $f$ is
well-defined on $M^{k}\times (S^{m-1}-S^{m-3})$ and also on $M^k\times
S^{m-3}$, if we ignore the undefined longitude coordinate. Observe that on each sphere
$*\times (S^{m-1}, S^{m-3})$, the map is just the rotation by angle $\tau(x)$
of the longitude coordinate. 
 
Roughly speaking, we are removing a bundle of trivial knots over
$M$ and replacing it with a bundle of non-trivial knots using a longitudinal twist determined by $\tau$. 

We now compute the knot groups of a frame twist-spun knot. If we spin the $m$-knot $K$ with singular set $\Sigma$, complement $C$, and boundary complement $X$, we represent the corresponding spaces for the spun knot $\sigma^{\psi,\tau}_M(K)$ by  $\sigma(\Sigma)$, $\sigma(C)$, and $\sigma(X)$.
The singular set $\sigma(\Sigma)$ of the spun knot  is just $M\times \Sigma$, and clearly the boundary knot complement $\sigma(X)$ of $\sigma^{\psi,\tau}_M(K)$ is simply $M\times X$. Therefore, if $G\cong \pi_1(X)$ is the boundary knot group for $K$, then $\pi_1(M)\times G$ is the boundary knot group for $\sigma^{\psi,\tau}_M(K)$. Of course if $\Sigma$ is not connected, then these calculations hold separately for each component.  
 
To calculate the group $\pi_1(\sigma(C))$ of the complement $\sigma(C)$ of the knot $\sigma^{\psi,\tau}_M(K)$, we observe as in \cite{GBF} that $\sigma(C)\cong Y\cup Z$, where \begin{align*}
Y&= S^n-(S^{n-2}\cup\text{int} (M^k\times D^m))\\
Z&= M^k\times (D_+^m-K_+)
\end{align*}
and $Y\cap Z= M^k\times (S^{m-1}-S^{m-3})$. But then $Y\sim_{h.e.} S^1$, $Z\sim_{h.e.} M\times C$, and $Y\cap Z\sim_{h.e.} M\times S^1$ (see \cite{GBF} or \cite{GBF1} for more details). Therefore, by the van Kampen theorem, $\pi_1(\sigma(C))\cong \Z *_{\pi_1(M)\times \Z} (\pi_1(M)\times \bG)$, where $\bG\cong \pi_1(C)$. 

We can simplify this expression by computing the homomorphisms $\pi_1(M)\times \Z\to \Z$ and $\pi_1(M)\times \Z\to \pi_1(M)\times \bG$ induced by the inclusions $Y\cap Z\into Y$ and $Y\cap Z\into Z$. To make things simpler, let us choose a basepoint which is on the boundary component $S^{m-1}-S^{m-3}$ of the knotted disk complement $D_+^m-K_+$ over some arbitrary base point of $M$. Then a simple meridian in  the trivial knot complement $S^{m-1}-S^{m-3}$ is also a simple meridian of   the spun knot $\sigma(K)$. Furthermore, by the construction this meridian represents  generators of  the  $\Z$  factor of $\pi_1(Y\cap Z)$ and of $\pi_1(Y)\cong \Z$. Also, it represents a simple meridian of $*\times K$, meaning the copy of $K$ attached over the basepoint of $M$ (note that $D_+^m-K_+$ is homotopy equivalent to $C$). Therefore, in the expression 
$$\pi_1(\sigma(C))\cong \Z *_{\pi_1(M)\times \Z} \pi_1(M)\times \bG,$$
the $\Z$ term on the left can be identified on the right with the product of some representative of some simple meridian in $\bG$ with the identity element of $\pi_1(M)$. 

Now, the map $\pi_1(M)\to  \pi_1(M)\times \bG$ given by the restriction to the first component of the map induced by the inclusion $Y\cap Z\to Z$ is simply the inclusion $\alpha\to \alpha\times 1$, while the map $\pi_1(M)$ to $\pi_1(Y)\cong \Z$ is determined by the map $\tau$. As computed in \cite{GBF} and \cite{GBF1}, if $x\in \pi_1(M)$, then the map $\pi_1(M)\to \Z$ determined  by the inclusion $M\times *\subset M\times (S^{m-1}-S^{m-3})\into Y$, where $*$ is our basepoint in $S^{m-1}-S^{m-3}$, is given by $\alpha\to \text{deg}(\tau(\alpha))$. But we have already observed that in $\pi_1(\sigma(C))$, the generator of $\Z\cong \pi_1(Y)$ can be identified with a certain simple meridian $\bg$ of $K$ and $\sigma(K)$; with the proper choice of meridian representing the element of weight one in $\bG$, this is the same element $\bg$ as before. Putting these computations together, we see that 
$$\pi_1(\sigma(C))\cong \frac{\pi_1(M)\times \bG}{<\{x^{-1}\bg^{\text{deg}(\tau(x))}\}>},$$
where the ``denominator'' is the normal subgroup generated by all $x^{-1}\bg^{\text{deg}(\tau(x))}$ as $x$ ranges over the generators of $\pi_1(M)$. Note, in particular, that if $\tau$ is the trivial map, then $\pi_1(\sigma(C))\cong \bG$, a well-know result for frame-spun and superspun knots (see \cite{Su92} and \cite{C70}).

Lastly, let us compute the map induced by inclusion from the boundary knot group $\pi_1(\sigma(X))$ to the knot group $\pi_1(\sigma(C))$. Assume that $\Sigma$ and $M$ are connected, and let us pick a  useful base point so that we fix the groups involved within their isomorphism classes. Clearly this base point must be within $\sigma(X)$, but in our calculation of $\pi_1(\sigma(C))$, we assumed a basepoint in $Y\cap Z\cong M\times (S^{m-1}-S^{m-3})$, where the second factor is the complement of the trivial boundary knot left after removing a trivial open disk knot $(D^m_{-}, D^{m-2}_{-})$ from the knot $K$ to be spun. We can let our choice of base in $M$ remain arbitrary. As for the second term, recall that in the spinning construction, we are free to remove from $K$ any trivial disk knot in a neighborhood of any locally-flat point. In particular, clearly we are free to choose $(D^m_{-},D^{m-2}_{-})$ and the neighborhood $\bar N(\Sigma)$ of the  singular set of $K$ so that $D^m_{-}\cap \bar N(\Sigma)=\emptyset$ and $(\bar D^m_{-},\bar D^{m-2}_{-})\cap \bar N(\Sigma)\cong (D^{m-1},D^{m-3}) $, in other words so that the boundary unknot  $(S^{m-1}, S^{m-3})$ intersects $(\bd \bar N(\Sigma), \bd \bar N(\Sigma)\cap K)$ in a trivial disk knot, say one of its hemispheres. With these choices, we can find base points that lie in both $M\times  (S^{m-1}-S^{m-3})$ and in $\sigma(X)$. It then follows immediately from the above calculations that the map induced by inclusion from $\pi_1(\sigma(X))$ to $\pi_1(\sigma(C))$ is the following composition:
\begin{equation*}
\begin{CD}
\pi_1(M)\times G@>\text{id}\times \phi >> \pi_1(M)\times \bG@>\text{projection}>> \frac{\pi_1(M)\times \bG}{<\{x^{-1}\bg^{\text{deg}(\tau(x))}\}>}.
\end{CD}
\end{equation*}
If $\Sigma$  has more than one component, then this formula holds for each component where $\phi$ is induced by the inclusion of the boundary knot corresponding to the appropriate component (and we choose base points appropriately). 

Let us formalize these calculations as a proposition:

\begin{proposition}\label{P: spin group}
Suppose $K$ is a not necessarily locally-flat knot $K^{m-2}\subset S^m$ with connected singular set $\Sigma$ and that $K$ is frame twist-spun  about a connected manifold $M^k$ embedded with framing in $S^{n-2}$ and with twisting function $\tau$. Suppose that $K$ has knot group $\bG$ and boundary knot group $G$. Then the boundary knot group of $\sigma^{\psi,\tau}_M(K)$ is $\pi_1(M)\times G$ and the knot group is $\pi_1(\sigma(C))\cong \frac{\pi_1(M)\times \bG}{<\{x^{-1}\bg^{\text{deg}(\tau(x))}\}>}$ for some element $\bg\in\bG$ of weight one. Furthermore, if $M$ and $\Sigma$ are connected, the map from the boundary knot group to the knot group is given by the composition  $\pi_1(M)\times G\overset{\text{id}\times \phi}{\to} \pi_1(M)\times \bG \to \frac{\pi_1(M)\times \bG}{<\{x^{-1}\bg^{\text{deg}(\tau(x))}\}>}$, where $\phi$ is induced by inclusion. If $\Sigma$ is not connected, then the boundary knot about the component $M\times \Sigma_i$ of the singular set of $\sigma(K)$ has group $\pi_1(M)\times \pi_1(X_i)$, where $X_i$ is the corresponding component of the boundary knot of $K$, and the map $\pi_1(M)\times \pi_1(X_i)\to \pi_1(\sigma(C))$  induced on fundamental groups by inclusion is given by the composition $\pi_1(M)\times \pi_1(X_i)\overset{\text{id}\times i_*}{\to} \pi_1(M)\times \bG \to \frac{\pi_1(M)\times \bG}{<\{x^{-1}\bg^{\text{deg}(\tau(x))}\}>}$, where $i_*:\pi_1(X_i)\to \pi_i(C)=\bG$ is induced by inclusion. 
\end{proposition}

Using this, we can prove the following realization statement:
\begin{theorem}\label{T: man real}
Given groups $G$ and $\bG$ which satisfy conditions \eqref{I: finite}-\eqref{I: closure} of Theorem \ref{T: class}, a connected manifold $M^k$ that  embeds with framing in $S^{n-2}$, $m=n-k\geq 6$, and a homomorphism $\eta:\pi_1(M)\to \Z$, there exists a knot $K\subset S^n$ with singular set $M$, boundary knot group $G\times \pi_1(M)$, knot group $\pi_1(M)\times \bG/<\{x^{-1}\bg^{\eta(x))}\}>$, and inclusion homomorphism between them given by the composition $\pi_1(M)\times G\overset{\text{id}\times \phi}{\to} \pi_1(M)\times \bG \overset{\text{proj.}}{\to} \frac{\pi_1(M)\times \bG}{<\{x^{-1}\bg^{\text{deg}(\eta(x))}\}>}$. 
\end{theorem}
\begin{proof}
By Theorem \ref{T: class} and Proposition \ref{P: nontriv}, we can construct an $m$-disk knot (with non-trivial boundary knot) for the given range of $m$ with groups $G$ and $\bG$ and map $\phi$. Coning on the boundary gives a sphere knot with point singularity and the same groups. We will frame twist-spin these knots about $M$ embedded in $S^{n-2}$ with some framing. To choose $\tau$, notice that the map $\pi_1(M)\to \Z$ gives a map in $\Hom(H_1(M), \Z)$ by factoring through the abelianization. This, in turn, gives an element of $H^1(M)$ and hence an element of $[M,K(\Z,1)]=[M,S^1]$. Any element in this homotopy class gives us a map $\tau$ which induces the appropriate map $\eta:\pi_1(M)\to \pi_1(S^1)\cong \Z$ (see \cite[Chptr. 8]{Sp}). Therefore, applying the fact that the isomorphism $\pi_1(S^1)\cong \Z$ is given by taking degrees of maps, the theorem follows from the previous proposition by frame twist-spinning if we can show that we are free to spin in such a way that the given element of weight one, $\bg\in \bG$, will play the role of  the simple meridian of the same label in the statement of the proposition. 

For this, notice that, in the proof of sufficiency in Theorem \ref{T: class}, the disk knot we end up with is the cocore of the $2$-handle attached along a curve representing $\bg$. Therefore, $\bg$ is a simple meridian of the disk knot. In particular, it bounds the core of the handle. So when, after coning the boundary, we split the knot into $D^n_+$ and $D^n_-$, let us take as $D^n_-$ a thin neighborhood of the core of the handle in the handle. Then $D^n_-$ is constructed around a locally flat point, as desired, and furthermore, it is clear that $\bg$ represents the meridian of the trivial knot $S^{n-3}\subset S^{n-1}$ in $\bd D^n_-$, at least up to choice of orientation. But this is the meridian which also represents the $\bg$ of Proposition \ref{P: spin group}. So, by making a proper choice of orientation for the generator of $\pi_1(Y)$ in the frame twist-spinning construction, we obtain a knot with the desired groups.
\end{proof}

Using these above constructions and Corollary \ref{C: s.c. man}, we obtain the following classification theorem:

\begin{comment} in a large range of dimensions the knot group pairs of PL knots whose singular sets consist of a single stratum that is a 
simply-connected manifold $\Sigma$ with $H_2(\Sigma)=0$ and such that $\Sigma$ can be embedded with framing in a sphere.
\end{comment}

\begin{theorem}\label{T: real sc sing}
If the PL sphere knot $K: S^{n-2}\into S^n$ can be stratified so that the singular set $\Sigma$ consists of a single stratum which is  a $2$-connected manifold, then
for the groups $\bG$ and $G$ to be the fundamental groups of the respective complements, $C$ and $X$,  of the knot and its boundary knot and for $\phi: G\to \bar G$ to be the homomorphism induced by inclusion, it is necessary and sufficient that $G$ and $\bar G$ satisfy the Kervaire conditions with elements of weight one $g\in G$ and $\bar g\in \bar G$ such that $\phi(g)=\bar g$. On the other hand, given such a manifold $\Sigma$  of dimension $m$ that can  be embedded with  framing into a sphere $S^{n-2}$, $n-2\geq m+4$, and $G$, $\bG$, and $\phi$ satisfying these conditions, there exists a PL knot $K: S^{n-2}\subset S^n$ whose singular set is $\Sigma$, whose groups are $G$ and $\bG$, and whose homomorphism induced by inclusion is $\phi$.
\begin{comment}
\begin{enumerate}
\item  $G$ and $\bG$ are finitely presentable,
\item  $G/[G,G]\cong \bG/[\bG,\bG]\cong \Z$,
\item  $H_2(G)=H_2(\bG)=0$,
\item  There exist elements  $g\in G$ and $\bar g \in \bG$ and a homomorphism $\phi: G\to \bG$ such that 
\begin{enumerate}
\item $G$ is the normal closure of $g$,
\item $\bG$ is the normal closure of $\bar g$,
\item $\phi(g)=\bg$.
\end{enumerate}
\end{enumerate}
\end{comment}
\end{theorem}
\begin{proof}
The necessity of the conditions on the boundary group $G$ is the content of Corollary \ref{C: s.c. man}, while that for $G$ has been shown in Theorem \ref{T: nec con}. The condition on the map follows from the remark in the proof of Corollary \ref{C: s.c. man} according to which we can take $g$ to be represented by a simple meridian of the boundary knot. This is also a simple meridian of $K$ under inclusion and hence its normal closure is all of $\bG$ as in the proof of Theorem \ref{T: nec con}.

For sufficiency,  we an apply Theorem \ref{T: man real}, noting that the triviality of $\pi_1(M)$ reduces the whole knot group and the boundary knot group to $\bG$ and $G$, respectively. 
\end{proof}

\begin{remark}
These constructions and computations can be taken slightly further by taking knot sums of the knots constructed in the theorem with each other, with locally-flat knots, or even with any other singular knots we might have available. Then we obtain knot groups which are the free products amalgamated along  appropriate meridians. 

\begin{comment}Of course in this case the boundary knots might have several components, but the group of each boundary knot component  simply maps into the appropriate term of the amalgamated free product by the map induced by inclusion (all modulo  appropriate choices of basepoints, of course). \end{comment}
\end{remark}

\begin{remark}
As a further alternate construction, we could frame-twist spin about a manifold $M$ that is not connected; we can even spin different knots about different components or use components of different dimensions. 
\begin{comment} If the knots being spun are locally-flat, then this is no different from taking the knot sum of separate spinnings, for once we have spun about one component of $M$,  the other manifold components can be ambiently isotoped within the knot into an arbitrarily small unknotted region. In other words, we could consider the intermediate knot (after just one spinning) as a knot sum with the trivial knot and then move any remaining components of $M$ into the trivial knot component by isotopy; then we can perform the next spin in this component to obtain a knot sum of spun knots, which is then clearly equivalent to the knot obtained by spinning about all components at once and then performing the isotopy which ``slides'' one part of the knot across the other. However, this will not in general be true if we spin knots that are not locally-flat. For example, if we spin a singular knot about two manifolds which are non-trivially linked, then in general we will not be able to slide the knotted pieces across each other since we can not slide nicely through a singular set of the embedding. Nonetheless, \end{comment}
By performing calculations similar to those above and inducting one component at a time, we can see that the knot groups will still be amalgamated free products of the form 
\begin{equation*}
\frac{\pi(M_1)\times \bG_1}{<\{x^{-1}\bg_1^{\text{deg}(\tau_1(x))}\}>}*_{r_1}\cdots*_{r_{m-1}} 
\frac{\pi(M_m)\times \bG_m}{<\{x^{-1}\bg_m^{\text{deg}(\tau_m(x))}\}>},
\end{equation*}
where the subscripts indicate components of $M$ and the relations $r_i$  equate the elements of each component of the product which correspond to the central meridian given in our above computations by the generator of  $\pi_1(Y)\cong \pi_1(S^1)$. \begin{comment}(To fix a common basepoint for each space so that the above formula makes sense, we may choose a basepoint in $Y$, link it to each component $Z_i$ of the spin by appropriate simple disjoint paths, and then compute homotopy groups after replacing each $Z_i$ with the union of $Z_i$ with the path, since clearly nothing has changed up to homotopy equivalence and hence up to isomorphisms of the fundamental groups.) \end{comment}
The technical details of this computation are left to the reader. 
\end{remark}

Lastly, as one more example of the kinds of groups we can realize,  
we compute the groups of the suspension of a knot. As usual, we let $(\bG, G)$ stand for the knot group/boundary knot group pair of the knot $K$. If we let $\mf{S}$ stand for suspensions (unfortunately, the symbols $\sigma$, $\Sigma$, and $S$ have already been employed), the suspension of a knot pair $(S^n, K)$ is the suspended pair $\mf{S}(K)=(\mf{S}S^n, \mf{S}K)$. Let us abuse notation and refer to the complement and boundary complement of  $\mf{S}(K)$ as $\mf{S}(C)$ and $\mf{S}(X)$, although these will not actually be suspensions. Since the suspension points are part of the knot, $ \pi_1(\mf{S}(C))\cong\pi_1(C)\cong \bG$. On the other hand, $\mf{S}(X)$ will be $C_+\cup (X\times I)\cup C_-$, two copies of $C$ attached to $X\times I$, one at $X\times 0$ and the other at  $X\times 1$, by the canonical inclusions. This is homotopy equivalent to two copies of $C$ identified along $X$. So, if $\Sigma $ is connected, $X$ will be connected and by the van Kampen theorem,  $\pi_1(\mf{S}(X))\cong \bG*_{G}\bG$, where the identification maps are those induced by inclusions $X\into C$. In other words, if $i_{\pm}: X_{\pm}\into C_{\pm}$ are two copies of the inclusion map and $i_{\pm*}:\pi_1(X_{\pm})\to \pi_1(C_{\pm})$ are the induced maps, then $\pi_1(\mf{S}(X))\cong \frac{\bG*\bG}{<i_{+*}(\gamma)[i_{-*}(\gamma)]^{-1}>}$ as the denominator runs over all $\gamma\in G$. From the geometry, the map $\pi_1(\mf{S}(X))\to \pi_1(\mf{S}(C))$ induced by inclusion is simply the projection $ \bG*_{G}\bG\to \bG$ induced by the identity on each factor.

\begin{comment}
Lastly, as one more example of the kinds of groups we can realize,  
we provide the groups of the suspension of a knot with a connected singular set. It follows from a straightforward computation that if $G$ and $\bar G$ are the knot groups of our original knot, then the group of the suspended knot will again be $\bar G$, while the boundary group of the suspended knot will be $\bar G*_G \bar G$. The homomorphism induced by inclusion is simply the projection induced by the identity on each factor. 
\end{comment}

\bibliographystyle{amsplain}
\bibliography{bib}

\end{document}